\documentclass[article,12pt]{elsarticle}
\makeatletter
\def\ps@pprintTitle{%
 \let\@oddhead\@empty
 \let\@evenhead\@empty
 \def\@oddfoot{\centerline{\thepage}}%
 \let\@evenfoot\@oddfoot}
\makeatother
\usepackage{epsfig}
\usepackage{amssymb,amsthm,amsmath,amsfonts,array,amscd}
\usepackage{geometry}
\usepackage{graphics}
\journal{Statistics and Probability Letters}

\begin{document}

\begin{frontmatter}

\title{\large LOCAL EFFICIENCY OF INTEGRATED GOODNESS-OF-FIT
TESTS UNDER SKEW ALTERNATIVES}

\author[ad1]{A. Durio}
\author[ad2,ad3]{Ya.Yu.Nikitin\corref{mycorrespondingauthor}}
\cortext[mycorrespondingauthor]{Corresponding author}

\ead{yanikit47@mail.ru}

\address[ad1]{Department of Economics "S. Cognetti de Martiis", University of Turin, Lungo Dora Siena 100/A, 10153, Torino, Italy}
\address[ad2]{Department of Mathematics and Mechanics, Saint Petersburg State University, 7/9 Universitetskaya nab., St. Petersburg, 199034 Russia}
\address[ad3]{National Research University - Higher School of Economics, Souza Pechatnikov, 16, St.Peters\-burg 190008, Russia}

\begin{abstract}
 The efficiency of distribution-free {\it integrated} goodness-of-fit
tests  was studied by Henze and Nikitin (2000, 2002) under location alternatives.
We calculate local Bahadur efficiencies of these tests
under more realistic generalized skew alternatives. They turn out to be
unexpectedly high.
\end{abstract}

\begin{keyword}
Integrated empirical process \sep Bahadur efficiency \sep skew alternative

 \MSC 62G10  \sep 62G20 \sep 62G30

\end{keyword}

\end{frontmatter}



\section{Introduction}

Goodness-of-fit testing is one of  the most important problems in
Statistics. If the hypothetical distribution is continuous, one
can apply distribution-free tests based on
functionals of the empirical process. Most known tests of such
type are the Kol\-mogorov and Cram\'er-von Mises tests and
their variants, see, e.g., \cite{ShoWel} and \cite{Niki}.

In search of new distribution-free tests with possibly better
efficiency properties, Henze and Nikitin \cite{HenNik00}, \cite{HenNik02} proposed new
test statistics based on the {\it integrated} empirical process.
They found their limiting distributions and calculated local Baha\-dur
efficiencies for location alternatives. These efficiencies
are comparable with the efficiencies of usual
distribution-free tests, but there exist also some interesting distinctions
in favor of these new tests. Gradually statistical inference using integrated
empirical processes becomes quite popular, see, e.g., \cite {Alv}, \cite{Bouz},
\cite{Jing} and \cite{Kuriki}.

However, the location alternative is a simplest alternative which
is not very realistic in practice, particularly because it preserves
the symmetry of the underlying distribution. In many situations it is more reasonable to assume 
asymmetric alternative models. The most interesting and simple example of such alternative models
in the case of normal distribution was introduced in \cite{Azza85}.
Let $\Phi$ and $\varphi$ denote the distribution
function and the density of the standard normal law.
Azzalini \cite{Azza85} proposed the skew-normal distribution
depending on the real parameter $\theta$ and having the density
$$
g(x,\theta) = 2 \varphi(x) \Phi(\theta x) \,,\,  x \in {\mathbb R}\, ,\,
\theta \geq 0.
$$
It is evident that for any $\theta$ the function $g(x,\theta)$ is
a density and that for $ \theta = 0$ we get the standard normal
density. Later the properties of Azzalini's skew-normal model and
its generalizations were considered in numerous papers.
Finally they were described and collected in \cite{AzzaCap}.

For any symmetric distribution function $F$ with the density $f$ and any symmetric distribution function $G$ with the density $g$
we can consider the {\it generalized skew distribution}
with the density
\begin{equation}
\label{altern}
h(x, \theta) = 2 f(x) G(\theta x)\,, \, x \in {\mathbb R}, \, \, \theta \geq 0.
\end{equation}

Note that this model is more general than that considered in \cite{DN02} and \cite{DN03}
in view of the emergence of almost arbitrary distribution function $G$ instead of initial distribution function $F$. This model is 
described and advocated in \cite{AzzaCap}.

It is quite interesting to calculate the efficiencies of
integrated distri\-bution-free tests mentioned above under
the generalized skew alternative (\ref{altern}). We select the Bahadur efficiency
as it is well-adapted for such calculations while other types of
efficiencies such as Pitman, Chernoff or Hodges-Lehmann are not applicable
or do not discriminate between two-sided tests. See \cite{Niki} for
details concerning the calculation of efficiencies and their
interrelations.

The calculation of local Bahadur efficiency of common
distribution-free tests under skew alternatives was performed in
\cite{DN02} and \cite{DN03}. In the present paper
we  calculate the efficiencies of the {\it integrated} tests
under the more general alternative (\ref{altern}).

 General expressions for local Bahadur efficiencies in case of one-parameter families of
alternatives can be found in \cite{Niki}. However we
cannot apply them as the alternative
(\ref{altern}) requires some additional analysis. This analysis was partially done in \cite{DN02}, \cite{DN03}.
We use corresponding results in sections 2 and 3 when calculating
the efficiencies for five examples of symmetric distributions with
different tail behaviors. These efficiencies are taken together in Table 1 of Section 4.
They demonstrate that the efficiencies of integrated tests are appreciably higher than of usual tests.
Section 5 is devoted to the analysis of local optimality of tests under consideration.

\section{Tests Based on Integrated Empirical Process.}

Let $X_1, ... , X_n$ be a random sample from the density $h(x,\theta)$
given by (\ref{altern}) and depending on the known symmetric
density $f$ and symmetric distribution function $G$, and a real parameter $\theta \geq 0$.
Let
\begin{equation}
\label{altern2}
H(x, \theta) = 2\int _{-\infty}^x f(u) G(\theta u) du\,, \, x \in {\mathbb R}, \, \, \theta \geq 0,
\end{equation}
be the distribution function corresponding to this density. We want to
test the goodness-of-fit hypothesis  $ H_0: \theta = 0 $ against
the alternative $H_1: \theta > 0$. Let $F_n $ be the empirical distribution function
based on the sample $X_1,...,X_n.$

Some well-known goodness-of-fit tests are based on the Kolmogorov
stati\-stic
$$
D_n = \sqrt{n} \sup_{t} |F_n (t) - F(t)|,
$$
on the Chapman -- Moses statistic
$$
\omega_n^1 = \sqrt{n} \int _{{\mathbb R}} (F_n(t)- F(t)) dF(t),
$$
on the Cram\'er -- von Mises statistic
$$
\omega_n^2 = n \int _{{\mathbb R}} \left(F_n (t) - F(t)\right)^2 dF(t),
$$
and on the Watson statistic
$$
U_n^2 = n \int _{{\mathbb R}} \left(F_n (t) - F(t)- \int _{{\mathbb R}} (F_n (s)- F(s))
dF(s)\right)^2 dF(t).
$$

These statistics are distribution--free and
can be considered as functionals of the empirical
processes
$$
\beta_n(x)= \sqrt{n} (F_n(x)-F(x))\,, \quad x \in {\mathbb R}  \,,
$$
or
$$
\alpha_n(u)= \sqrt{n} (G_n(u)- u))\,, \quad 0 \leq u \leq 1\,,
$$
where the empirical distribution function $ G_n$ is based on the uniform sample
$F(X_j), j=1,...,n.$  Clearly $\beta_n(x)=\alpha_n(F(x))$, and we
can write
\begin{multline*}
D_n =  \sup_{x} |\beta_n(x)|= \sup_{u} |\alpha_n(u)|, \qquad  \omega_n^1 = \int _{{\mathbb R}} (\beta_n(x)) dF(x)= \int_0^1 \alpha_n(u) du,\\
\omega_n^2 = \int _{{\mathbb R}} (\beta_n(x))^2 dF(x)= \int_0^1 \alpha_n^2(u) du\,,\\
U_n^2 = \int_{{\mathbb R}} \left(\beta_n(x) -\int_{{\mathbb R}}\beta_n(s)dF(s) \right)^2
dF(x)= \int_0^1 ( \alpha_n(u) -  \int_0^1 \alpha_n(s)ds)^2 du\,.
\end{multline*}

Henze and Nikitin, see \cite{HenNik00} and \cite{HenNik02},  proposed similar
but more complicated statistics based
on the {\it integrated empirical process} and studied their Bahadur
local efficiency for the location alternative. Let
$$
\bar{F}_n(x)= \int_{-\infty}^x F_n(t) dF(t), \qquad
\bar{F}(x)=\int_{-\infty}^x F(t) dF(t)= \frac{1}{2} F^2(x)
$$
denote the {\it integrated } empirical distribution function and the {\it integrated }
hypothetical distribution function respectively. Then the  integrated empirical
process is
$$
B_n(x)= \sqrt{n} [\bar{F}_n(x)-\bar{F}(x)]= \int_{-\infty}^x
\beta_n(t) dF(t), \quad  x \in {\mathbb R},
$$
while the integrated uniform empirical process becomes
$$
A_n(u)=  \int_{0}^u \alpha_n(s) ds\,, \qquad 0 \leq u \leq 1\,.
$$
The integrated analogs of the classical statistics $D_n$,
$\omega_n^1,$ $\omega_n^2$ and $U_n^2$ were defined in \cite{HenNik00, HenNik02} as
$$
\begin{array}{ll}
\vspace{10pt}
 \bar{D}_n &= \sup_{x} |B_n(x)|= \sup_{u} |A_n(u)|\,,\\
 \vspace{10pt}
\bar{\omega}_n^1 &=\int_{{\mathbb R}} B_n(t) dF(t)= \int_0^1 A_n(u)
du,\quad \bar{\omega}_n^2 = \int_{{\mathbb R}} B_n^2(t) dF(t)=
\int_0^1 A_n^2(u) du,\\
 \bar{U}_n^2 &=\int_{{\mathbb R}} \left(B_n(t) - \int_{{\mathbb R}}B_n(s)dF(s)\right)^2
dF(t)= \int_0^1 (A_n(u) - \int_0^1 A_n(s)ds)^2 du.
\vspace{10pt}
\end{array}
$$

Henze and Nikitin in \cite{HenNik00} and \cite{HenNik02} derived limiting distributions,
large deviation asymptotics, local Bahadur efficiencies for location
alternatives, and studied the conditions of local Bahadur optimality for these
statistics. In next sections we will carry through this program under the generalized
skew alternative (\ref{altern}).

\section{Bahadur local efficiency: general expressions}

In the rest of the paper, we consider alternative (1) with the symmetric density $f$ having finite variance. The distribution function $G$ and the density $g = G'$ are assumed to be symmetric as well. They all satisfy the following conditions.

\smallskip

\noindent{\it Condition 1}. We require that the density $g$ with $g(0) >0$ is positive and differentiable within its support.
By symmetry we always have $g'(0)=0$.

\noindent{\it Condition 2}. Let $f$ and $g$ be such that uniformly in $x
\in {\mathbb R}$
$$
 H(x,\theta) -F(x) \sim 2\theta  g(0)\int_{-\infty}^{x}uf(u)du, \quad as \quad \theta \to 0,
$$
where $\sim$ is the usual sign of equivalence.

\noindent{\it Condition 3}. Suppose that
 $$
 K(\theta) \sim 2
g^2(0)\int_{{\mathbb R}}x^2 f(x)dx \ \theta^2 ,\, \, as \quad \theta
\to 0,
$$
 where $ K(\theta)$ is the well-known
Kullback -- Leibler information \cite{Bah71}
$$ K(\theta) := \int_{{\mathbb R}} \ln
\{ h(x,\theta)/h(x,0) \}h(x,\theta) dx= 2\int_{{\mathbb R}} \ln \{
2G(\theta x)\} f(x)G(\theta x)dx.
$$
These conditions are very natural and are valid for various densities $f$ and $g.$
Condition 2 was obtained by using the Taylor expansion of $G(\theta x)$ for small $\theta$
and extracting the leading term. To get the Condition 3, we use the expansion 
$$ 
y \ln y =  y-1  -\frac12 (y-1)^2 +o\{ (y-1)^2\}, \mbox{as} \, \ y \to 1,
$$
which implies as $\theta \to 0$, for any $x$  (since $g'(0)=0$)
$$
2 G (\theta x) \ln \{ 2 G(\theta x)\}=  2 G(\theta x) -1 +\frac12 \{ 2 G(\theta x)-1\}^2 +o(\theta^2) = 2g(0)\theta x +2g^2(0)\theta^2 x^2 +o(\theta^2).
$$

Substituting this in the definition of $K(\theta)$ above and integrating, we get under weak additional requirements the Condition 3.

It is not difficult to impose sufficient conditions on $f$ and $g$
ensuring such behavior but we prefer the formulation of regularity conditions
in form of Conditions 1-3.

Now we describe in short the definition and calculation of Bahadur efficiency.
Details can be found in \cite{Bah67}, \cite{Bah71}, and \cite{Niki}.

Suppose that $T = \{ T_n \}$ is a sequence of statistics,  such
that as $ n \to \infty$
 $$
 \begin{array}{ll}
 \vspace{6pt}
 &a) \quad T_n   \longrightarrow  b(T,\theta) \qquad \mbox{in probability
under} \, \  H_1; \\
&b) \quad  n^{-1} \ln P( T_n \geq \varepsilon) \longrightarrow   -r( T,\varepsilon) \quad \mbox{under} \, \ H_0,
\end{array}
$$
where the function $r(T, \varepsilon)$ is continuous in $\varepsilon$ for sufficiently small $ \varepsilon >0.$ Condition a) is a variant of the law of large numbers under $H_1$
while condition b) is always non-trivial and describes the (logarithmic) large deviation
behavior of test statistics under the null-hypothesis. Then the
exact Bahadur slope is defined as
$$
  c(T,\theta) = 2 r(T, b(T,\theta))\,,
$$
while the local Bahadur efficiency is defined by
$$
 e^{B} (T) = \lim_{\theta \to 0+ } \frac{
c(T,\theta)}{2K(\theta)}\, .
$$

In all the examples considered in this paper we have
\begin{equation}
\label{index}
c(T,\theta) \sim l(T, f)  4g^2(0)  \theta^2 \,,\quad \mbox{as} \quad \theta
\to 0+,
\end{equation}
 where the functional $ l(T, f) $ is called the
local index. Then we have
 \begin{equation}
 \label{calc}
 e^{B} (T) = \frac{ l(T,f)}{\sigma^2(f)},
\end{equation}
where $\sigma^2(f)$ is the variance of the density $f.$

For our test statistics $\bar{D}_n,\, \bar{\omega}_n^1,\,
\bar{\omega}_n^2$ and $\bar{U}_n^2 $ the function $b(T,\theta)$
was found in \cite{HenNik00} and \cite{HenNik02} in terms of alternative
distribution function $H(x,\theta):$
$$
b(\bar{D},\theta) \equiv \sup_s |\int_{-\infty}^s
(H(x,\theta)-F(x)) dF(x)|\,,
$$

$$
b(\bar{\omega}^1,\theta) \equiv \int_{{\mathbb R}} \left[
\int_{-\infty}^s ( H(x,\theta)-F(x)) dF(x)\right] dF(s)\,,
$$

$$
b(\bar{\omega}^2,\theta) \equiv \int_{{\mathbb R}} \left[
\int_{-\infty}^s (H(x,\theta)-F(x)) dF(x)\right]^2 dF(s)\,,
$$
$$
b(\bar{U}^2,\theta) \equiv \int_{{\mathbb R}} \left[ \int_{-\infty}^s
(H(x,\theta)-F(x)) dF(x)\right]^2 dF(s)-
$$
$$
-\left(\int_{{\mathbb R}} \left[ \int_{-\infty}^s ( H(x,\theta)-F(x))
dF(x)\right] dF(s)\right)^2\,.
$$

\medskip

Using (2), regularity conditions 1 - 3, and setting
\begin{equation}
\label{nota}
v(x) = \int_{-\infty}^{x}uf(u)du,\, \quad  q(s)=\int_{-\infty}^{s} v(x) f(x) dx,
\end{equation}
we easily arrive to the
following expressions for the local representations of functions
$b$ as $\theta \to 0+:$
$$
\begin{array}{ll}
\vspace{6pt}
&b(\bar{D},\theta) \sim 2 \theta g(0) \sup_s |q(s)|,\quad b(\bar{\omega}^1,\theta) \sim  2 \theta g(0) \int_{{\mathbb R}} q(s) f(s)ds,\\
\vspace{6pt}
&b(\bar{\omega}^2,\theta) \sim 4 \theta^2 g^2 (0) \int_{{\mathbb R}} q^2 (s)f(s) ds, \\
&b(\bar{U}^2,\theta) \sim 4 \theta^2 g^2(0)\left[ \int_{{\mathbb R}}
q^2 (s)f(s) ds - \left(\int_{{\mathbb R}}  q(s) f(s) ds\right)^2\right]\,.
\end{array}
$$

Applying the large deviation asymptotics of integrated statistics
from \cite{HenNik00} and \cite{HenNik02}, we find the following local behavior of exact slopes for our
 test statistics as $\theta \rightarrow 0+:$
$$
\begin{array}{ll}
\vspace{6pt}
&c(\bar{D}, \theta) \sim 12 b^2 (\bar{D},\theta), \quad c(\bar{\omega}^1, \theta) \sim 45 b^2 (\bar{\omega}^1,\theta)\,,\\
\vspace{6pt}
&c(\bar{\omega}^2, \theta) \sim \mu_0 b(\bar{\omega}^2,\theta) \, \  \mbox{with}\,  \mu_0=
31.2852..., \quad  c(\bar{U}^2, \theta) \sim \pi^4 b(\bar{U}^2,\theta).
\end{array}
$$

Combining these formulas with the asymptotics of
functions $b$ given above, we easily obtain the expressions for
the local exact indices $l(T,f)$, see (\ref{index}), of our statistics. The factor $4g^2(0)$ disappears
when calculating the local efficiency according to (\ref{calc}). Hence we may write
\begin{equation}
\label{main}
e^B (T) = \frac{l(T,f)}{\sigma^2(f)}.
\end{equation}
We get now the following expressions for local indices of our statistics:
$$
\begin{array}{ll}
&l(\bar{D},f)= 12 \sup_s q^2(s)\,, \quad l(\bar{\omega}^1,f)= 45 \left(\int_{\mathbb R} q(s) f(s) ds \right)^2, \,
l(\bar{\omega^2},f)= \mu_0 \int_{{\mathbb R}}  q^2(s) f(s) ds,\\
&l(\bar{U}^2,f)= \pi^4\left(
\int_{{\mathbb R}} q^2(s) f(s) ds- \left(\int_{{\mathbb R}} q(s) f(s) ds\right)^2 \right)^2.
\end{array}
$$
Note that the efficiencies {\do not depend on $G.$}

\section{Bahadur local efficiency: examples and discussion}

We will  calculate local indices for following five standard symmetric densities $f:$
$$
\begin{array}{lll}
&f_1(x)= (2 \pi)^{-1/2} \exp (-x^2/2), \qquad \qquad &( \mbox{normal density})\\
&f_2(x)= e^x/(1+e^x)^2, \qquad \qquad &( \mbox{logistic density})\\
&f_3(x)= 1/(\pi(1-x^2)^{1/2}){\bf 1}_{[-1,1]} (x)\,,\,  \qquad \qquad &( \mbox{arcsine density})\\
&f_4(x)= \frac{1}{2} {\bf 1}_{[-1,1]} (x), \qquad \qquad &( \mbox{uniform density})\\
&f_5(x)= 8/(3\pi (1+x^2)^3), \qquad \qquad &(\mbox{non-standardized  Student-5 density.}) \\
\end{array}
$$
Using the notation (\ref{nota}) for all $f_i, i=1,...,5,$ we see that
$$
\begin{array}{lll}
\vspace {10pt} v_1(x) &= -\frac{1}{\sqrt{2\pi} } e^{-x^2/2}, x \in {\mathbb R}, \qquad &v_2(x) = -\ln(1+ e^{x})+\frac{x e^{x}}{1+e^{x}}, x
\in {\mathbb R}, \\
\vspace {10pt} v_3(x) &= -\frac{1}{\pi}\sqrt{1-x^2},\, -1
\leq x \leq 1, \qquad &v_4(x) = -\frac{1}{4}( 1-x^2 ), \,-1\leq
x \leq 1, \\
\vspace {10pt} v_5(x) &= -\frac{2}{3\pi(1+x^2)^2}, x \in {\mathbb R}.
\end{array}
$$
Next we calculate for our densities the functions $q_i, \ i=1,\dots,5:$
$$
\begin{array}{lll}
\vspace {10pt} q_1(s) &=  - \frac{\Phi(s\sqrt{2})}{2\sqrt{\pi}}\,,s \in {\mathbb R}, \quad q_2(s) =\frac{1+e^s+se^{2s}- (e^{2s}-1)\ln(1+e^s)}{2(1+e^s)^2}-\frac{1}{2},\,s \in {\mathbb R}, \\
\vspace {10pt} q_3(s) &=-\frac{s+1}{\pi^2},\, |s|\leq 1, \,  \qquad q_4(s) =\frac{s^3-3s-2}{24},\ |s|\leq 1,\\
\vspace {10pt} q_5(s) &=- \frac{s(279+511s^2+385s^4+105s^6)+105(1+s^2)^4 \arctan(s)}{216\pi^2(1+s^2)^4}-\frac{35}{144 \pi},\, s \in {\mathbb R}.
\end{array}
$$

Now we proceed to the calculation of local indices for our five densities.
Obser\-ving that $\sup_s |q_i(s)|$ are respectively $1/(3 \pi)$,
1/2, $2/\pi^2$, $1/6$ and $35/(72\pi),$ we obtain
$$
l(\bar{D}_n,f_1)= 0.95493, \,\, l(\bar{D}_n,f_2)=3, \,\,
l(\bar{D}_n,f_3)=48/\pi^4, $$
$$ l(\bar{D}_n,f_4)=1/3, \,\,
l(\bar{D}_n,f_5)=1225/(432
\pi^2).
$$

Since $\int_{-\infty}^{+\infty}q_i(s)f_i(s) ds$, for $1\le i \le 5$, are respectively
$1/(4 \sqrt{\pi})$, $-1/4$, $-1/\pi^2$, $-1/12$ and $-35/(144 \pi)$
we obtain
$$
l(\bar{\omega}_n^1,f_1)= 0.8952, \,\,
l(\bar{\omega}_n^1,f_2)=45/16,\,
\, l(\bar{\omega}_n^1,f_3)=45/\pi^4,
$$
$$
l(\bar{\omega}_n^1,f_4)=5/16, \, \,
l(\bar{\omega}_n^1,f_5)=6125/(2304 \pi^2).
$$

Finally knowing that $\int_{-\infty}^{+\infty} q_i^2(s) f_i(s)
ds$ are respectively $0.02914$, 0.09107, $3/(2
\pi^4)$, $13/1260$ and $1225/(62208 \pi^2)+ (46189+39200
\pi^2 )/(663552
\pi^4),$ we obtain
$$
l(\bar{\omega}_n^2,f_1)= 0.91154, \,
\, l(\bar{\omega}_n^2,f_2)= 2.84924, \, \, l(\bar{\omega}_n^2,f_3)=
0.48176,
$$
$$
l(\bar{\omega}_n^2,f_4)= 0.32278, \,
\, l(\bar{\omega}_n^2,f_5)= 0.27204.
$$

According to (\ref{main}) we need also the variances $\sigma^2(f)$ which are in our cases respec\-ti\-vely 1, $\pi^2 /3$, 1/2, 1/3 and 1/3. We  summarize our calculations in Table 1 where for comparison we also report the local efficiency of classical statistics $D_n$, $\omega_n^1,$ $\omega_n^2$ and $U_n^2$ given in \cite{DN03} for skew alternatives corresponding to the same five densities.

The inspection of
 this table and its comparison with  Table 3 in \cite[p.80]{Niki}
 and corresponding tables in \cite{HenNik00} and \cite{HenNik02} shows that the ordering of tests is similar
to the location case. This is favorable for practitioners: they
seldom know the structure of the alternative but can use the same
test both for the location and skew models.

\begin{table}
\renewcommand{\arraystretch}{1.2}
\begin{center}
\caption{ \normalsize Local Bahadur efficiencies under skew
alternatives. }
\vspace{6pt}
\begin{tabular}{|c|c|c|c|c|c|c|c|} \hline
Statistic &\multicolumn{5}{|c|}{Distribution}\\ \hline & Gauss &
Logistic & Arcsine & Uniform & Student-5\\ & & & & &   \\
 \hline $D_n$
& 0.637 & 0.584 &0.810 & 0.750 & 0.540\\
 \hline $\omega^1_n$ & 0.955 & 0.912 &
0.985 & 1 & 0.862\\
 \hline $\omega^2_n$ & 0.907 & 0.855 & 1 &
0.987 & 0.802 \\
  \hline $U_n^2$& 0.486 & 0.420 & 0.662 & 0.658 &
0.373  \\
 \hline $\bar{D}_n$
& 0.955 & 0.912 &0.985 & 1 & 0.862 \\
 \hline $\bar{\omega}_n^1$
& 0.895 & 0.855 & 0.924& 0.938 & 0.808\\
 \hline $\bar{\omega}_n^2$
& 0.912 & 0.866 & 0.963 & 0.968 & 0.816\\
 \hline $\bar{U}_n^2$& 0.900 & 0.846  & 1 & 0.986 & 0.792 \\
 \hline
\end{tabular}
\end{center}
\end{table}

However the efficiencies of integrated statistics are in
most cases {\it considerably higher than of classical ones.} This justifies
the use of integrated statistics for skew alternatives.

Note that the efficiencies of the statistics $\bar{D}_n$ and
$\omega_n^1$ coincide. It is not surprising as they have the same
local indices. It explains the maximal efficiency 1
attained by  $\bar{D}_n$ for the uniform
distribution, while for  $\omega_n^1$ the same
 was discovered in \cite{DN03}.  Another curious observation is that
for the normal law the  efficiencies under location
and skew alternatives coincide. This is a characteristic
 property of the normal law, see \cite{DN03}.
The efficiency 1 for $\bar{U}_n^2$ for the arcsine density is unexpected and will be interpreted
below.

Note that the so-called Pitman limiting relative efficiency of the
considered statistics is equal to the local Bahadur efficiency
under somewhat stronger regularity conditions.
It can be verified in the same way as in \cite{Wie} and \cite{Niki}.

Lachal in an interesting paper \cite{Lachal} studied $p$-fold integrated empirical processes and corresponding statistics. He considered, however, only location alternatives. For $p=0$ his results coincide with the conclusions of \cite{HenNik00} and \cite{HenNik02}. Moreover, for $p>1$ his tests demonstrate the decrease of efficiency (found numerically) when $p$ grows, but the theoretical calculations are hardly possible.

\section{Conditions of local optimality.}

As is well known \cite{Bah67}, \cite[Ch.6]{Niki} the
local asymptotic optimality (LAO) of a sequence $\{T_n\}$ in
Bahadur sense means that $e^B(T)=1$ or, by (\ref{calc}), one has
\begin{equation}
\label{cond}
 l(T,f) = \int _{\mathbb R} x^2 f(x)dx.
 \end{equation}

We are interested  in those densities $f$ when (\ref{cond}) is
true; such densities under corresponding regularity conditions form the so-called {\it domain of LAO}.
The study of this "inverse" problem was  started by Nikitin (1984). The {\it a priori} regularity conditions are described in \cite[Ch.6]{Niki}, 
we underline the assumption $f(x) >0$  for all $x.$  In the sequel $C_1, C_2,  \dots$ denote some indefinite non-null real constants.

Note first of all that $\rho(s) := \int_{-\infty}^{s} \int_{-\infty}^{x}uf(u)du f(x) dx $ attains its maximum for $s=\infty.$ Indeed, the extremum condition is $\rho'(s) = f(s) \int_{-\infty}^{s}uf(u)du =0,$ and as $f>0,$ we see that $\rho'(s) =0 $ only for $s=\infty.$

Let apply this argument for the Kolmogorov statistic. Due to symmetry of $f$, we get, integrating by parts and applying the Cauchy-Schwarz inequality, that
\begin{multline*}
l(D,f)=  12\sup_s \left(\int_{-\infty}^{s} v(x) f(x) dx\right)^2 = 12\sup_s \left(\int_{-\infty}^{s} \int_{-\infty}^{x}uf(u)du f(x) dx\right)^2 = \\= 12\left(\int_{\mathbb R} \int_{-\infty}^{x}uf(u)du  f(x) dx\right)^2 = 12 \left(\int_{\mathbb R} u (F(u) -\frac12) f(u) du \right)^2 \leq \\
 \leq  12 \int_{-\infty}^{\infty}u^2 f(u) du
 \int_{-\infty}^{\infty}(F(u)- \frac12)^2 dF(u) = \int _{\mathbb R} x^2 f(x)dx.
\end{multline*}
Hence the condition of LAO (\ref{cond}) in virtue of the condition of equality in Cauchy-Schwarz inequality reduces
to the condition
\begin{equation}
\label{eq}
 F(x) - 1/2 = C_{1} x
\end{equation}
on the support of $f$. This implies that $f$ is constant on a
symmetric interval around zero. We
consider this as a characterization of the symmetric
uniform distribution.

We remark that the local optimality of the same statistic $\bar{D}_n$ under
the  location alternative is valid  for logistic distribution,
see \cite{HenNik00}, this emphasizes the difference between
these two types of alternatives.

The arguments for the sequence $\{\bar{\omega}_n^1\}$ are
similar but the result is different. We have, using integration by parts, the symmetry of the density $f$ and the Cauchy-Schwarz inequality
\begin{multline*}
\quad l(\bar{\omega}^1,f)= 45 \left(\int_{\mathbb R} q(s) f(s) ds \right)^2 = 45 \left(\int_{\mathbb R}\int_{-\infty}^{s} v(x)  f(x) dx\ f(s)ds\right)^2 =\\ = 45 \left(\int_{\mathbb R} v(x) (1- F(x)) f(x) dx\right)^2 =\frac{45}{4}\left( \int_{\mathbb R} v(x) d((1-F(x))^2 \right)^2=\\=
\frac{45}{4}\left( \int_{\mathbb R}x \left((1-F(x))^2 -\frac13\right) f(x) dx\right)^2 \leq \\ \leq  \frac{45}{4} \int_{0}^1 (z^2 -1/3)^2 dz
\int_{\mathbb R}x^2 f(x) dx = \sigma^2(f).
\end{multline*}

Using the condition of equality in Cauchy-Schwarz inequality, we see that the condition of LAO  is valid iff
\begin{equation}
\label{equa}
(1-F(x))^2 -\frac13  = C_2 x
\end{equation}
on the support of symmetric $f.$ This is impossible, unlike (\ref{eq}), since for symmetric distribution function $F$ we have $F(0)=\frac12$, and this contradicts the equation (\ref{equa}).

For the integrated statistic $\bar{\omega}_n^2$ such direct arguments are problematic. Therefore we will apply the general theory developed in \cite[Ch.6]{Niki}. According to it, any sequence of statistics $\{T_n\}$ defines the "leading function" $v_{T}$ (or sometimes a set of them) which specifies the most efficient direction in the space of alternatives $H(x,\theta).$  To describe the domain of LAO we need to solve the equation
$$
H'_{\theta} (x,0) = C_3 v_T (F(x)) \, \mbox{with some constant } \, C_3.
$$
The set of alternatives $H(x,\theta)$ should satisfy some regularity conditions listed and discussed in \cite[Ch.6]{Niki}. The skew family
(\ref{altern2}) under conditions 1-3 satisfies them for a very broad set of densities $f$ and distribution function's $G$. Hence we can apply this theory subject to knowledge of "leading functions" which can be at times very involved. For the integrated statistic $\bar{\omega}_n^2$ the set of leading functions was found in \cite{HenNik00} by variational methods and consists of eigenfunctions of some boundary-value problem, namely
$$
\psi_j(x) \ = \ \cos \kappa_j \ \sinh \left( \kappa_j (1-x)) + \cosh (\kappa_j) \sin (\kappa_j (1-x)\right),\ j\ge 1,
$$
with $\kappa_j$ being the consecutive positive zeros of the equation $\tan(x) + \tanh(x) =0.$
Consider the first of these functions $\psi_1.$ It does not change its sign on $[0,1].$ Hence the distribution function $F$ of interest for us has to satisfy the
differential equation
$$
\int_{-\infty}^x u f(u) du \ = C_4 \left( \cos \kappa_1 \ \sinh \left( \kappa_1 (1-F(x)\right) + \cosh
 \kappa_1 \sin (\kappa_1 (1-F(x))\right).
$$
Differentiating this equation, we can obtain on the support of $f$ an implicit equation for $F$ but we are not able to obtain its explicit solution.

It is curious that the more complicated integrated statistic $\bar{U}_n^2$ has a much simpler domain of LAO. The leading functions here
 \cite{HenNik02} are $\sin(\pi jx), j=1,2,...$ Only the first function keeps the sign on $[0,1]$ so that we arrive to the differential equation
$$
\int_{-\infty}^x u f(u) du = C_5 \sin \pi F(x), \, x \in {\mathbb R}.
$$

After differentiation we get the equation
$$
f(x) ( x - C_6\cos \pi F(x)) =0,
$$
which results on the set $\{x: f(x) \neq 0\}$ in the solution
$$
F(x) = 1 - \pi^{-1} \arccos(x/C_7) = \pi^{-1}\arcsin(x/C_6) + 1/2, \, -C_7 \leq x \leq C_7,
$$
corresponding to the symmetric arcsine density 
$$ f(x) =  \left(\pi \sqrt{C_7^2 -x^2}\right)^{-1} {\bf 1}\{ -C_7 \leq x \leq C_7\}.$$
It may be observed that we got a characterization of arcsine density by the property of LAO for $\bar{U}_n^2$
  under the skew alternative. This explains the appearance of 1 in the last row in Table 1 above.

\section{Acknowledgement}

The research of second author was supported by RFBR grant No. 16-01-00258.

\end{document}